\documentclass[10pt]{article}

\usepackage{amsmath}
\usepackage{amsfonts}   
\usepackage{amssymb}  

\newlength{\originalbase}
\newtheorem{theorem}{Theorem}[section]
\newtheorem{proposition}[theorem]{Proposition}

\newtheorem{corollary}[theorem]{Corollary}

\input amssym.def
 
\input amssym.tex 

\renewcommand{\SS}{{\Bbb S}}
\newcommand{\RR}{{\Bbb R}}

\newcommand{\HH}{{\Bbb H}}

\newcommand{\XX}{{\Bbb X}}

\newcommand{\qed}{\vrule height1.8ex width0.9ex depth0.0ex}

\newcommand{\euler}{\raisebox{.4ex}{$\chi$}}

\begin{document}

\input epsf

\pagestyle{plain}
\pagenumbering{arabic}

\begin{center}
{\large\bf Bonnesen-type inequalities for surfaces of constant curvature}\\
\end{center}

\begin{center}
{\bf Daniel A. Klain}\footnote{Research supported in part by NSF grant
\#DMS-9803571.} \\
Department of Mathematical Sciences\\
University of Massachusetts Lowell\\
Lowell, MA 01854 USA  \\
dklain@cs.uml.edu\\
\end{center}
\vspace{0.6 cm}
\begin{quotation}
{\small
\noindent
{\bf Abstract \/}   A {\em Bonnesen-type} inequality is a
sharp isoperimetric
inequality that includes an error estimate in terms of 
inscribed and circumscribed regions.
A kinematic technique is used to
prove a Bonnesen-type inequality for  
the Euclidean sphere (having constant
Gauss curvature $\kappa > 0$) 
and the hyperbolic plane (having constant 
Gauss curvature $\kappa < 0$).  These generalized
inequalities each converge to the classical Bonnesen-type inequality
for the Euclidean plane as $\kappa \rightarrow 0$.
}
\end{quotation}

\vspace{1 cm}
\noindent
{\large\bf Introduction}

A {\em Bonnesen-type} inequality is a
sharp isoperimetric
inequality that includes an error estimate in terms of 
inscribed and circumscribed regions.  The classical example
runs as follows:
\begin{quotation}
Suppose that $K$ is a compact convex set in $\RR^2$.
Denote by $A_K$ and $P_K$ the area and perimeter of $K$
respectively.
Let $R_K$ denote the circumradius
of $K$, and let $r_K$ denote the inradius of $K$.  Then
\end{quotation}
\begin{equation}
P_K^2 - 4\pi A_K \geq \pi^2 (R_K - r_K)^2.
\label{bonneu}
\end{equation}

\noindent
The classical isoperimetric inequality immediately follows, namely,
\begin{equation}
P_K^2 - 4\pi A_K \geq 0,
\label{isoper}
\end{equation}
with equality if and only if $R_K = r_K$, that is, if and only if $K$ is a Euclidean disc.
Proofs of these inequalities, along with variations and
generalizations, can be found in any of 
\cite{BonnIneq,Santa,Osser}, for example.

In this note the
kinematic methods of Santal\'o and Hadwiger are used to
prove Bonnesen-type inequalities for  
the Euclidean sphere (having constant
Gauss curvature $\kappa > 0$) 
and the hyperbolic plane (having constant 
Gauss curvature $\kappa < 0$).  Section 1 outlines necessary 
background material from integral geometry.
In Section 2 we derive the first of the two main theorems
in this article, a Bonnesen-type inequality for the sphere,
stated in Theorem~\ref{bonns}.
The second main theorem of this article, Theorem~\ref{bonn}, is a
Bonnesen-type inequality for the hyperbolic plane, 
derived in Section 3.  The limiting case as $\kappa \rightarrow 0$ in either of
Theorems~\ref{bonns} and~\ref{bonnh} yields 
the classical Bonnesen inequality~(\ref{bonneu}), as
described above.
A brief and direct proof of~(\ref{bonneu}) 
using kinematic arguments, also described in \cite{Santa}, 
is presented at the close of Section 1
as a contrast to those of the subsequent sections.

\section{\large\bf Background: Integral geometry of surfaces}
Denote by $\XX_{\kappa}$ the surface of constant curvature $\kappa$,
specifically:
$$\XX_{\kappa} = \left\{
\begin{array}{ccl}
\hbox{Euclidean 2-sphere of radius } 1/\sqrt{\kappa} & \hbox{if} & \kappa > 0 \\
\hbox{Euclidean plane }\RR^2 & \hbox{if} & \kappa = 0 \\
\hbox{Hyperbolic plane of constant curvature } \kappa & \hbox{if} & \kappa < 0 \\
\end{array}\right.$$

A compact set $P \subseteq \XX_{\kappa}$ is a {\em convex polygon} if $P$
can be expressed a finite intersection of closed half-planes (or closed hemi-spheres
in the case of $\kappa > 0$, with the added requirement that 
$P$ lie in inside an open hemisphere).  A {\em polygon} is a finite union of convex polygons.
More generally, a set $K \subseteq \XX_{\kappa}$ will be called {\em convex} 
if any two points of $K$ can be connected by a line segment
inside $K$, where the notion of line segment is again suitably
defined for each context (spherical, Euclidean, hyperbolic).  
For $\kappa >0$, a convex set is again required to lie inside an open hemisphere.
Denote by ${\cal K}(\XX_{\kappa})$ the set of all {\em compact} convex sets in
$\XX_{\kappa}$.  

For $K \in {\cal K}(\XX_{\kappa})$, denote by $A_K$ the {\em area} of $K$, and
denote by $P_K$ the {\em perimeter} of $K$.  If $\dim K = 1$ then $P_K$
is equal to {\em twice} the length of $K$.  (This assures that perimeter $P$
is continuous in the Hausdorff topology on compact sets in $\XX_{\kappa}$.)

If $K$ is a finite
union of compact convex sets in $\XX_{\kappa}$, denote by 
$\euler_{K}$ the {\em Euler
characteristic} of $K$.  If $K$ is a compact convex set, then $\euler_{K} = 1$
whenever $K$ is nonempty, while $\euler_{\emptyset} = 0$.  More generally, $\euler$
extends to all finite unions of compact convex sets via iteration of the
inclusion-exclusion identity:
$$\euler_{K \cup L} + \euler_{K \cap L} = \euler_{K} + \euler_{L}.$$

Our primary tool for studying inequalities will be the {\em principal kinematic
formula} \cite[p. 321]{Santa} for compact convex sets in $\XX_{\kappa}$.

\begin{theorem}[Principal Kinematic Formula for $\XX_{\kappa}$]
For all finite unions $K$ and $L$ of compact convex sets in
$\XX_{\kappa}$,
\begin{equation}
\int_{g} \euler_{K \cap gL} \; dg \; = \; 
 \euler_{K} A_L + \frac{1}{2\pi} P_K P_L
+ A_K \euler_{L} - \frac{\kappa}{2\pi} A_K A_L.
\label{hpkfeq}
\end{equation}
\label{hpkf}
\end{theorem}
The integral on the left-hand side of~(\ref{hpkfeq})
is taken with respect to area on
$\XX_{\kappa}$ and
the invariant Haar
probability measure on the group $G_0$ of isometries of $\XX_{\kappa}$
which fix a base point $x_0 \in \XX_{\kappa}$.  To 
define this more precisely, denote by $t_x$ the unique 
translation of $\XX_{\kappa}$ 
(or minimal rotation, in the case of $\kappa > 0$)
that maps $x_0$ to a point $x \in \XX_{\kappa}$.
Then define
\begin{equation}
\int_{g} \euler_{K \cap gL} \; dg = 
\int_{x \in \XX_{\kappa}} \int_{\gamma \in G_0} 
\euler_{K \cap t_x(\gamma L)} \; d\gamma \; dx,
\label{kinint}
\end{equation}
where we use the probabilistic normalization
$$\int_{\gamma \in G_0} d\gamma = 1.$$

The classical proof of Theorem~\ref{hpkf} can be found in \cite{Santa}.
For a valuation-based proof of Theorem~\ref{hpkf}, see  
\cite{Lincee} (for the Euclidean and spherical cases) and
\cite{Klain-Hyperbolic} (for the hyperbolic plane).
Surveys and other recent work on kinematic formulas in convex, integral, 
and Riemannian geometry and
their applications include \cite{Fu,Howard,Lincee,Santa,Sch-Wie,Zhang-dkf}.

Choose a fixed base point $x_0 \in \XX_{\kappa}$. 
For $r \geq 0$, denote by $D_r$ the set of points in $\XX_{\kappa}$
that lie at most a distance $r$ from $x_0$.  We will refer to 
$D_r$ as the {\em disc of radius r} in $\XX_{\kappa}$. 

Recall that, for $\kappa \neq 0$,
\begin{equation}
P_{D_r} = \frac{2\pi}{\sqrt{\kappa}} \sin (\sqrt{\kappa}r) 
\;\;\; \hbox{ and } \;\;\; 
A_{D_r} = \frac{2\pi}{\kappa}(1 - \cos(\sqrt{\kappa}r)).
\label{diskforms}
\end{equation}
See, for example, \cite[p. 85]{Still}.  The limiting
cases as $\kappa \rightarrow 0$
yield the Euclidean formulas $P_{D_r} = 2\pi r$ and 
$A_{D_r} = \pi r^2$.

Theorem~\ref{hpkf} leads in turn to the following version of
Hadwiger's containment theorem for convex subsets of surfaces 
\cite{Had-Contain1,Had-Contain2,Lincee,Santa}. 

\begin{theorem}[Hadwiger's Containment Theorem]
Let $K,L \in {\cal K}(\XX_{\kappa})$ with non-empty interiors.  
If
\begin{equation}
P_K P_L \leq 2\pi (A_K + A_L) - \kappa A_K A_L,
\label{cineq}
\end{equation}
then there exists an isometry $g$ of $\XX_{\kappa}$
such that either $gK \subseteq L$ or 
$gL \subseteq K$. 
\label{contain}
\end{theorem}

\noindent
{\bf Proof:\/}  First, consider the case in which $K$ and $L$
are convex {\em polygons} in $\XX_{\kappa}$.
Suppose that, for every isometry $g$,
we have $gK \nsubseteq \hbox{int}(L)$
and $gL \nsubseteq \hbox{int}(K)$.  In this instance,
whenever the $K$ and $gL$ overlap, the boundary intersection
$\partial K$ and $\partial L$ will consist of a discrete set
of 2 or more points (except for a measure zero set of motions $g$).
In other words, for almost all isometries $g$,
$$\euler_{\partial K \cap g\partial L} \geq 2\euler_{K \cap gL}.$$
On integrating both sides with respect to $g$, it follows from
the kinematic formula~(\ref{hpkfeq}) that
$$\frac{1}{2\pi} P_{\partial K} P_{\partial L}
\geq 2 \left(
A_K + \frac{1}{2\pi}P_KP_L + A_L - \frac{\kappa}{2\pi} A_K A_L
\right).$$
Recall that $P_{\partial K} = 2P_K$, 
and similarly $P_{\partial L} = 2P_L$, so that
$$\frac{1}{\pi} P_K P_L
\geq  
A_K + \frac{1}{2\pi}P_KP_L + A_L - \frac{\kappa}{2\pi} A_K A_L.$$
Hence,
$$P_KP_L \geq 2\pi( A_K + A_L ) - \kappa A_K A_L.$$
In other words, if~(\ref{cineq}) holds with {\em strict} inequality ($<$)
then there exists an isometry $g$ such that
either $gK \subseteq \hbox{int}(L)$
or $gL \subseteq \hbox{int}(K)$.  Since the set
$\{(K,L,g) \; | \; gK \subseteq L \hbox { or } gL \subseteq K \}$
is closed in the Hausdorff topology, the theorem also holds for
the case of equality in~(\ref{cineq}), as well as for
compact convex sets $K$ and $L$ in $\XX_{\kappa}$ that are not polygons.
\qed

We will also make use of the following elementary fact about the perimeter of compact convex
subsets of $\XX_{\kappa}$.

\begin{proposition}  Suppose that $K, L \in {\cal K}(\XX_{\kappa})$ and suppose that $K \subseteq L$.
Then $P_K \leq P_L$.
\label{monoper}
\end{proposition}
In other words, perimeter is monotonic on compact convex sets in $\XX_{\kappa}$.

Evidently Proposition~\ref{monoper} is not true for arbitrary (non-convex) sets.  
Nor does it hold for convex-like subsets of the sphere that 
do not lie inside an open hemipshere.

\noindent
{\bf Proof:\/}  Suppose that $\kappa \leq 0$ (so that we consider the Euclidean or hyperbolic plane).
Let $N$ be a line segment of length $d$, so that $P_N = 2d$.  
Note that $\euler_{N \cap K} = 1$ if and only if $N \cap K \neq \emptyset$; otherwise
$\euler_{N \cap K} = 0$.  Moreover, $\euler_{N \cap K} \leq \euler_{N \cap L}$ since $K \subseteq L$.
On averaging over all motions of $N$, the kinematic formula, Theorem~\ref{hpkf}, implies that
$$A_K + \frac{1}{2\pi}P_N P_K \leq A_L + \frac{1}{2\pi}P_N P_L,$$
so that
$$\frac{A_K}{2d} + \frac{1}{2\pi} P_K \leq \frac{A_L}{2d} + \frac{1}{2\pi} P_L,$$
for all $d > 0$.  Taking the limit as $d \rightarrow \infty$ yields $P_K \leq P_L$.

For $\kappa > 0$ (the Euclidean sphere) replace the line segment $N$ 
with a great circle $C$.  
Recall that convex sets in the sphere are required to lie inside an open hemisphere,
so that 
$\euler_{C \cap K} = 1$ if and only if $C \cap K \neq \emptyset$, and
$\euler_{C \cap K} \leq \euler_{C \cap L}$ whenever $K \subseteq L$.
Because $\euler_{C} = A_C = 0$, 
the kinematic formula of Theorem~\ref{hpkf} implies that
$$\frac{1}{2\pi}P_C P_K \leq \frac{1}{2\pi}P_C P_L,$$
so that $P_K \leq P_L$ once again.
$\qed$

After setting $\kappa = 0$ in Theorem~\ref{contain}, a kinematic proof of 
the classical Bonnesen-type inequality~(\ref{bonneu}) is straightforward.

\noindent
{\bf Proof of the inequality~(\ref{bonneu}):\/}
If $K$ is a disc, 
then both sides of the inequality~(\ref{bonneu}) are equal to zero.

Suppose that $K \in {\cal K}(\RR^2)$ is not a disc, 
so that $r_K < R_K$.  For $r_K < \epsilon < R_K$ we can apply Theorem~\ref{contain}
to $K$ and the disc $D_{\epsilon}$ to obtain
$$P_K P_{D_{\epsilon}} > 2\pi (A_K + A_{D_{\epsilon}}).$$
It follows that
$$P_K 2\pi\epsilon 
> 2\pi (A_K + \pi \epsilon^2).$$
In other words,
$$
f(\epsilon) = -\pi \epsilon^2 + \epsilon P_K - A_K > 0,
$$
for all $\epsilon \in (r_K, R_K)$.  Since the leading coefficient of the
quadratic polynomial $f(\epsilon)$ is negative, it follows that
$f$ has two distinct roots, separated by the interval $(r_K, R_K)$.

Hence,  $\delta(f)/\pi^2 \geq (R_K - r_K)^2$, 
where $\delta(f)$ is the discriminant of $f$.  In other words,
$$P_K^2 - 4\pi A_K \geq \pi^2(R_K - r_K)^2.$$
$\qed$

\section{\large\bf Isoperimetry in $\SS^2$}

In this section we consider the case $\XX_{\kappa}$ for $\kappa > 0$.
For simplicity of notation, we first consider the case $\kappa = 1$.
A restatement of the main results for
general constant curvature $\kappa > 0$ 
is then given at the end of the section;
the proofs are entirely analogous to the case $\kappa = 1$.

Denote by $\SS^2$ the Euclidean unit sphere in $\RR^2$.
For $K \in {\cal K}(\SS^2)$ define the {\em circumradius} $R_K$
to be the greatest lower bound of all radii $R$ such that
some spherical disc of radius $R$ contains $K$.  Similarly,
define the {\em inradius} $r_K$ to be the least upper bound 
of all radii $r$ such that $K$ contains a spherical disc 
(i.e. spherical cap) of radius $r$.
Evidently $r_K \leq R_K$, with equality if and only if $K$ is a spherical disc.
Our restriction that a convex set must always lie in an open hemisphere
implies that $R_K < \frac{\pi}{2}$.

We will use Theorem~\ref{contain} to prove the following 
Bonnesen-type inequality for the sphere $\SS^2$.
\begin{theorem}[Bonnesen-type inequality for $\SS^2$]
Suppose $K \in {\cal K}(\SS^2)$.   Then
\begin{equation}
P_K^2 - A_K (4\pi -  A_K) \geq 
\frac{(\sin R_K - \sin r_K)^2
((2\pi -  A_K)^2 +  P_K^2)^2}
{4 (2\pi -  A_K)^2}
\label{bonnineqs}
\end{equation}
\label{bonns}
\end{theorem}

The inequality~(\ref{bonnineqs}) has the following simplification that also
provides equality conditions.
\begin{corollary}[Simplified Bonnesen-type inequality for $\SS^2$]
Suppose $K \in {\cal K}(\SS^2)$.  Then
\begin{equation}
P_K^2 - A_K(4\pi -  A_K) \geq 
\frac{1}{4}(\sin R_K - \sin r_K)^2
(2\pi -  A_K)^2.
\label{bonnineqs2}
\end{equation}
with equality if and only if $K$ is a spherical disc.
\label{bonns2}
\end{corollary}

\noindent
{\bf Proof of Corollary~\ref{bonns2}:} Since $ P_K^2 \geq 0$,
\begin{equation}
\frac{(\sin R_K - \sin r_K)^2
((2\pi -  A_K)^2 +  P_K^2)^2}
{4 (2\pi -  A_K)^2} 
\geq \frac{1}{4}(\sin R_K - \sin r_K)^2
(2\pi -  A_K)^2.
\label{yeah}
\end{equation}
The inequality~(\ref{bonnineqs2}) now follows from~(\ref{bonnineqs})
and ~(\ref{yeah}).

Equality holds in~(\ref{bonnineqs2}) and~(\ref{yeah}) 
if and only if either $P_K = 0$, in which case $K$ is a single point,
or if $R_K = r_K$, in which case $K$ must be a spherical disc.
\qed

The right-hand sides of~(\ref{bonnineqs}) and~(\ref{bonnineqs2}) are always
non-negative and are equal to zero if and only if $R_K = r_K$, that is,
if and only if $K$ is a disc.
These observations yield the following classical result as an immediate corollary.
\begin{corollary}[Isoperimetric inequality for $\SS^2$]
For $K \in {\cal K}(\SS^2)$,
$$
P_K^2 \geq A_K \left( 4\pi - A_K \right),
$$
with equality if and only if $K$ is a spherical disc.
\label{isoper2s}
\end{corollary}

Note that the complement $K'$ of $K$ in $\SS^2$, while not convex according
to our definition, has the same boundary and perimeter as $K$, while the inradius
and circumradius exchange roles.   Meanwhile,
$A_{K'} + A_K = 4\pi$, so that Corollaries~\ref{bonns2} 
and~\ref{isoper2s} are transformed as follows.

\begin{corollary}[Alternate simplified Bonnesen-type inequality for $\SS^2$]
Suppose $K \in {\cal K}(\SS^2)$.    Then
\begin{equation}
P_K^2 -  A_K A_{K'} \geq 
\frac{1}{16} \, (\sin R_K - \sin r_K)^2 (A_K - A_{K'})^2,
\label{bonnineqs2a}
\end{equation}
so that, in particular,
$$P_K^2 - A_K A_{K'} \geq 0$$
with equality in both cases if and only if $K$ is a spherical disc.
\label{bonns2a}
\end{corollary}

\noindent
{\bf Proof of Corollary~\ref{bonns2a}:} Since $ A_{K'} = 4 \pi -  A_K$
the left-hand sides of~(\ref{bonnineqs2}) and~(\ref{bonnineqs2a}) are the same. 
Meanwhile,
\begin{eqnarray*}
(A_K - A_{K'})^2 &=&  A_K^2 +  A_{K'}^2 - 2 A_K A_{K'} \\
&=&  A_K^2 + (4\pi - A_K)^2 - 2 A_K (4\pi -  A_K) \\
&=& 4 A_K^2 - 16\pi  A_K + 16 \pi^2 \\
&=& 4(2\pi -  A_K)^2,
\end{eqnarray*}
so that the 
right-hand sides of~(\ref{bonnineqs2}) and~(\ref{bonnineqs2a}) are the same
as well.
$\qed$

We now prove the main inequality of this section, Theorem~\ref{bonns}.

\noindent
{\bf Proof of Theorem~\ref{bonns}:}  If $K$ is a disc, 
then $R_K = r_K$, and both sides of~(\ref{bonnineqs}) are equal to zero.

Suppose that $K \in {\cal K}(\SS^2)$ is not a disc, so that $r_K < R_K$.  
For $r_K < \epsilon < R_K$ we can apply Theorem~\ref{contain}
to $K$ and the disc $D_{\epsilon}$ to obtain
$$P_K P_{D_{\epsilon}} > 2\pi (A_K + A_{D_{\epsilon}}) -  A_K A_{D_{\epsilon}}$$
It follows from~(\ref{diskforms}) that
\begin{eqnarray*}
P_K \sin  \epsilon
&>& A_K + 2\pi(1- \cos \epsilon) - 
A_K (1- \cos \epsilon)\\
&=& 2\pi (1 - \cos \epsilon) + A_K\cos \epsilon.
\end{eqnarray*}
Setting $P = P_K$ and $A = A_K$, we have
\begin{equation} 
P \sin  \epsilon - 2\pi 
> (A - 2\pi) \cos  \epsilon.
\label{ap}
\end{equation}
In order for $K \subseteq \SS^2$ to be convex, $K$ must be contained
in a hemisphere, so that $P \leq 2\pi$. It follows that
$$ P \sin  \epsilon - 2\pi \leq 0,$$
so that both sides of~(\ref{ap}) are {\em non-positive.}
Set $x = \sin \epsilon$, so that 
$\cos \epsilon = \sqrt{1-x^2}$.
Squaring both sides of~(\ref{ap}) reverses the order, yielding
\begin{eqnarray*}
P^2 x^2 - 4\pi  P x + 4 \pi^2 &<& (A - 2\pi)^2(1-x^2) \\
&=& -(2\pi -A)^2 x^2 + (4 \pi^2 - 4\pi A + A^2),
\end{eqnarray*}
so that
\begin{equation}
f(x) = [(2\pi -  A)^2 +  P^2]x^2 - 
4\pi  P x + (4\pi -  A)  A < 0,
\label{qs}
\end{equation}
for all $x \in (\sin r_K, \sin R_K)$. 

Since $K$ is not a disc, $K$ is not a point, so $P_K > 0$.  It follows that
$(2\pi -  A)^2 +  P^2 \geq P^2 > 0$, so that
the quadratic polynomial $f(x)$ defined by~(\ref{qs}) 
has a positive leading coefficient, and $f(x) > 0$ for
sufficiently large $|x|$.  Since $f(x) < 0$ for
$x \in (\sin  r_K, \sin  R_K)$, 
it follows that $f(x)$ has two
{\em real} roots, and that these two roots must lie on 
{\em different sides}
of the open interval $(\sin r_K, \sin R_K)$.

The discriminant $\delta(f)$ of the quadratic polynomial $f$ is computed as follows:
\begin{eqnarray*}
\delta(f) &=& (4\pi P)^2  - 4((2\pi -  A)^2 +  P^2)(4\pi -  A) A \\
&=& 4 [4\pi^2 P^2 - ((2\pi -  A)^2 +  P^2)(4\pi- A) A]\\
&=& 4 [4\pi^2 P^2 - 4\pi  P^2 A +  P^2 A^2 - 
(2\pi- A)^2(4\pi- A) A]\\
&=& 4 [P^2 (4\pi^2  - 4\pi   A +  A^2) - 
(2\pi- A)^2(4\pi- A) A]\\
&=& 4 [P^2 (2\pi -  A)^2 - 
(2\pi- A)^2(4\pi- A) A]\\
&=& 4 (2 \pi -  A)^2 (P^2 - A(4\pi -  A))
\end{eqnarray*}
The squared distance between the roots of a quadratic polynomial $f$
is given by its discriminant $\delta(f)$ divided by the square
of its leading coefficient.  Hence,
\begin{equation}
\frac{4 (2\pi-  A)^2(P^2 - A(4\pi - A))}{((2\pi -  A)^2 + P^2)^2}
\geq (\sin R_K - \sin r_K)^2.
\label{bs}
\end{equation}
which implies the inequality~(\ref{bonnineqs}). $\qed$

For the general case, 
denote by $\SS^2_\kappa$ the Euclidean sphere having radius $\frac{1}{\sqrt\kappa}$
and Gauss curvature $\kappa$.  In this case
our restriction that a convex set must always lie in an open hemisphere
implies that $R_K < \frac{\pi}{2\sqrt{\kappa}}$.  Note also that if $K$ and
$K'$ are complements in $\SS^2_\kappa$ then $A_K + A_{K'} = \frac{4\pi}{\kappa}$.

The inequalities of
this section are now summarized in full generality.  These generalized versions
follow immediately from the theorems above via a scaling argument. Alternatively
these more general cases can be proved in direct analogy to the proof given 
above for the case $\kappa = 1$.

\begin{theorem}[Bonnesen-type inequalities for $\SS_\kappa^2$]
Suppose $K \in {\cal K}(\SS_\kappa^2)$,
and $K'$ denote the complement of $K$ in $\SS_\kappa^2$.
Then the following inequalities hold:
\begin{equation}
P_K^2 - A_K (4\pi - \kappa A_K) \geq 
\frac{(\sin \sqrt{\kappa}R_K - \sin \sqrt{\kappa}r_K)^2
((2\pi - \kappa A_K)^2 + \kappa P_K^2)^2}
{4\kappa (2\pi - \kappa A_K)^2}
\label{s1}
\end{equation}

\begin{equation}
P_K^2 - A_K(4\pi - \kappa A_K) \geq 
\frac{1}{4\kappa}(\sin \sqrt{\kappa}R_K - \sin \sqrt{\kappa}r_K)^2
(2\pi - \kappa A_K)^2
\label{s2}
\end{equation}

\begin{equation}
P_K^2 \geq A_K \left( 4\pi - \kappa A_K \right)
\label{s3}
\end{equation}

\begin{equation}
P_K^2 - \kappa A_K A_{K'} \geq 
\frac{\kappa}{16} \, (\sin \sqrt{\kappa}R_K - \sin \sqrt{\kappa}r_K)^2 (A_K - A_{K'})^2
\label{s4}
\end{equation}

Equality holds in~(\ref{s2}),~(\ref{s3}), and~(\ref{s4}) if and only if $K$ is a spherical disc.
\label{omnibus-s}
\end{theorem}
Note that as $\kappa \rightarrow 0^+$ the inequalities~(\ref{s1}) and~(\ref{s2}) of 
Theorem~\ref{omnibus-s} yield the classical Bonnesen-type inequality~(\ref{bonneu})
for the Euclidean plane, while~(\ref{s3}) reduces to the classical isoperimetric
inequality~(\ref{isoper}).

\section{\large\bf Isoperimetry in $\HH^2$}

In this section we consider the case of constant negative curvature; that is,
$\kappa < 0$.  To simplify notation, we first consider the case $\kappa = -1$.
A restatement of the main results for
general constant curvature $\kappa < 0$ 
is then given at the end of the section;
the proofs are entirely analogous to the case $\kappa = -1$.

Let $\HH^2$ denote the hyperbolic plane having constant negative curvature -1.
For $K \in {\cal K}(\HH^2)$ define the {\em circumradius} $R_K$
to be the greatest lower bound of all radii $R$ such that
some hyperbolic disc of radius $R$ contains $K$.  Similarly,
define the {\em inradius} $r_K$ to be the least upper bound 
of all radii $r$ such that $K$ contains a hyperbolic disc of radius $r$.
Evidently $r_K \leq R_K$, with equality if and only if $K$ is a hyperbolic disc.

The following theorem is a limited analogue of the spherical 
Bonnesen-type inequality of Theorem~\ref{bonns}.
\begin{theorem}
Suppose $K \in {\cal K}(\HH^2)$.  
If $(2\pi +  A_K)^2 -   P_K^2 \geq 0$, then
\begin{equation}
P_K^2 - A_K(4\pi +  A_K) \geq 
\frac{(\sinh R_K - \sinh r_K)^2
((2\pi +  A_K)^2 -  P_K^2)^2}{4 (2\pi +  A_K)^2}. 
\label{bonnineq}
\end{equation}
\label{bonn}
\end{theorem}
The proof of Theorem~\ref{bonn} is deferred to the end of this section.
Although~(\ref{bonnineq}) appears almost identical to the spherical 
Bonnesen inequality~(\ref{bonnineqs}), up to change of sign in a few places, 
on more careful examination some other important differences appear.

Note that the inequality~(\ref{bonnineq}) may fail to hold if 
$(2\pi +  A_K)^2 -   P_K^2 < 0$.  For example, if $K$
is a line segment of length $c$, then $A_K = r_K = 0$, while $P_K = 2c$ and
$R_K = c/2$.  In this instance, the left-hand side of~(\ref{bonnineq})
is $O(c^2)$, while the right-hand side of~(\ref{bonnineq}) grows
exponentially in $c$.  This apparent deficiency will be addressed 
by Theorem~\ref{bonnh}.
  
Meanwhile, note that
the condition $(2\pi +  A_K)^2 -   P_K^2 \geq 0$ is not as strange as it may appear, 
when compared carefully to the spherical case of Theorem~\ref{bonns}.  
In the sphere we required that
a convex set be contained in a hemisphere, that is, a spherical disc of radius 
$\frac{\pi}{2\sqrt{\kappa}}$.   According to~(\ref{diskforms}), this spherical disc satisfies
$$P_{D_\frac{\pi}{2\sqrt{\kappa}}} = 
\frac{2\pi}{\sqrt{\kappa}} \sin (\sqrt{\kappa}\frac{\pi}{2\sqrt{\kappa}}) 
= \frac{2\pi}{\sqrt{\kappa}} \sin \frac{\pi}{2} = \frac{2\pi}{\sqrt{\kappa}},$$
since $\sin \frac{\pi}{2} = 1$.
The next corollary involves an analogous assumption
for the hyperbolic plane.  In this context our
replacement for the value $\pi/2$ will be $\eta = \sinh^{-1}(1) = \ln(1 + \sqrt{2})$.

\begin{corollary}
Suppose $K \in {\cal K}(\HH^2)$. 
If $K \subseteq D_\eta$, where $\sinh \eta = 1$, then
$$
P_K^2 - A_K (4\pi +  A_K) \geq 
\frac{(\sinh R_K - \sinh r_K)^2
((2\pi +  A_K)^2 -  P_K^2)^2}{4 (2\pi +  A_K)^2}.  
$$
\label{bonn2}
\end{corollary}

\noindent
{\bf Proof of Corollary~\ref{bonn2}:}  Recall from Proposition~\ref{monoper} that
perimeter $P_K$ is monotonic with respect to set inclusion 
when applied to {\em convex} sets.
If $K \subseteq D_\eta$ it follows from~(\ref{diskforms}) that
$$P_K \leq P_{D_\eta} 
= \frac{2\pi}{i} \sin (i\eta)
= 2\pi \sinh \eta = 2\pi,$$
where $i^2 = -1$.  It follows that
$$(2\pi +  A_K)^2 -   P_K^2 \geq 
4\pi^2 + 4\pi A_K + A_K^2 - 4\pi^2
= 4\pi A_K + A_K^2 \geq 0,$$
so that Theorem~\ref{bonn} applies.
$\qed$

If, contrary to the hypothesis of Theorem~\ref{bonn}, we
have $(2\pi +  A_K)^2 -  P_K^2 < 0$, then
\begin{equation}
 P_K^2 -  A_K(4\pi +  A_K) >
(2\pi +  A_K)^2 -  A_K (4\pi +  A_K) = 4\pi^2.
\label{a}
\end{equation}
Combining~(\ref{a}) with Theorem~\ref{bonn} yields the following theorem.

\begin{theorem}
Suppose $K \in {\cal K}(\HH^2)$.  
Then
\begin{equation}
P_K^2 - A_K(4\pi +  A_K) \geq \min \left( 4\pi^2,
\frac{(\sinh R_K - \sinh r_K)^2
((2\pi +  A_K)^2 -  P_K^2)^2}{4 (2\pi +  A_K)^2}  
\right)
\label{bonnineq2}
\end{equation}
\label{bonnh}
\end{theorem}

The right hand side of~(\ref{bonnineq2}) is always
non-negative and is equal to zero if and only if $R_K = r_K$, that is,
if and only if $K$ is a disc.
These observations yield the following corollary.
\begin{corollary}[Isoperimetric inequality for $\HH^2$]
For $K \in {\cal K}(\HH^2)$,
$$
P_K^2 \geq A_K \left( 4\pi +  A_K \right).
$$
Equality holds if and only if $K$ is a hyperbolic disc.
\label{isoper2}
\end{corollary}

\noindent
{\bf Proof of Theorem~\ref{bonn}:}  If $K$ is a disc, then $R_K = r_K$ 
and both sides of~(\ref{bonnineq}) are equal to zero.

Suppose that $K \in {\cal K}(\HH^2)$ is not a disc, 
so that $r_K < R_K$. For $r_K < \epsilon < R_K$ we can apply Theorem~\ref{contain}
to $K$ and the disc $D_{\epsilon}$ to obtain
$$P_K P_{D_{\epsilon}} > 2\pi (A_K + A_{D_{\epsilon}}) +  A_K A_{D_{\epsilon}}.$$
It follows from~(\ref{diskforms}) that
$$P_K \sinh  \epsilon 
> A_K + 2\pi (\cosh \epsilon - 1) + 
A_K (\cosh \epsilon - 1),$$
so that
\begin{equation}
 P_K \sinh  \epsilon + 2\pi 
> (2\pi +  A_K) \cosh  \epsilon.
\label{um}
\end{equation}
Set $x = \sinh \epsilon$, so that 
$\cosh \epsilon = \sqrt{1+x^2}$.
To simplify the notation, let $P = P_K$ and $A = A_K$.  Since
the right hand side of~(\ref{um}) is positive, we can square
both sides of~(\ref{um}) to obtain
$$ P^2 x^2 + 4\pi  P x + 4 \pi^2 > (2\pi +  A)^2(1+x^2)
= (2\pi +  A)^2 x^2 + (4 \pi^2 + 4\pi  A +  A^2),$$
so that
\begin{equation}
f(x) = [ P^2 - (2\pi +  A)^2]x^2 + 
4\pi  P x - (4\pi +  A)  A > 0,
\label{q}
\end{equation}
for all $x \in (\sinh  r_K, \sinh  R_K)$. 

Recall that, by the hypothesis of the theorem,
$ P^2 \leq (2\pi +  A)^2$.  

If $ P^2 < (2\pi +  A)^2$, then
the quadratic polynomial $f(x)$ defined by~(\ref{q}) 
has a negative leading coefficient, so that $f(x) < 0$ for
sufficiently large $|x|$.  Since $f(x) > 0$ for
$x \in (\sinh  r_K, \sinh  R_K)$, 
it follows that $f(x)$ has two
{\em real} roots, and that these two roots must lie on 
{\em different sides}
of the open interval $(\sinh r_K, \sinh R_K)$.

The discriminant $\delta(f)$ of $f$ is computed as follows:
\begin{eqnarray*}
\delta(f) &=& (4\pi P)^2  + 4( P^2 - (2\pi +  A)^2)(4\pi +  A) A \\
&=& 4 [4\pi^2 P^2 + ( P^2 - (2\pi+ A)^2)(4\pi+ A) A]\\
&=& 4 [4\pi^2 P^2 + 4\pi  P^2 A +  P^2 A^2 - 
(2\pi+ A)^2(4\pi+ A) A]\\
&=& 4 [P^2 (4\pi^2  + 4\pi   A +   A^2) - 
(2\pi+ A)^2(4\pi+ A) A]\\
&=& 4 [P^2 (2\pi +  A)^2 - 
(2\pi+ A)^2(4\pi+ A) A]\\
&=& 4  (2 \pi +  A)^2 (P^2 - A(4\pi +  A))
\end{eqnarray*}
The squared distance between the roots of a quadratic polynomial $f$
is given by its discriminant $\delta(f)$ divided by the square
of its leading coefficient.  Therefore,
$$
\frac{4  (2\pi +  A)^2(P^2 - A(4\pi +  A))}{((2\pi +  A)^2 -  P^2)^2}
\geq (\sinh R_K - \sinh r_K)^2,
$$
from which~(\ref{bonnineq}) then follows.

Finally, if $ P^2 = (2\pi +  A)^2$, then the right-hand 
side of~(\ref{bonnineq}) is zero,
while the left-hand side satisfies
$$P^2 - A(4\pi +  A) = 
\left(  P^2 -  A(4\pi +  A) \right) =
(2\pi +  A)^2 -  A(4\pi +  A) = 4\pi^2 > 0.
$$
\qed

For the general case of constant negative curvature $\kappa < 0$, 
let $\lambda = |\kappa|$, and let 
$\HH^2_\lambda$ denote the hyperbolic plane having Gauss curvature $-\lambda$.

The inequalities of
this section are now summarized in full generality.  These generalized versions
follow immediately from the theorems above via a scaling argument. Alternatively
these more general cases can be proved in direct analogy to the proof given 
above for the case $\kappa = -1$ (that is, $\lambda = 1$).

\begin{theorem}
Suppose $K \in {\cal K}(\HH^2_{\lambda})$.  
If $(2\pi + \lambda A_K)^2 -  \lambda P_K^2 \geq 0$, then
\begin{equation}
P_K^2 - A_K(4\pi + \lambda A_K) \geq 
\frac{(\sinh \sqrt{\lambda}R_K - \sinh \sqrt{\lambda}r_K)^2
((2\pi + \lambda A_K)^2 - \lambda P_K^2)^2}{4\lambda (2\pi + \lambda A_K)^2}. 
\label{h1}
\end{equation}
More generally, if $K \in {\cal K}(\HH^2_{\lambda})$ then
$$
P_K^2 - A_K(4\pi + \lambda A_K) \geq \min \left( \frac{4\pi^2}{\lambda},
\frac{(\sinh \sqrt{\lambda}R_K - \sinh \sqrt{\lambda}r_K)^2
((2\pi + \lambda A_K)^2 - \lambda P_K^2)^2}{4\lambda (2\pi + \lambda A_K)^2}  
\right)
$$
\label{omnibus-h}
\end{theorem}

In particular, if $K \subseteq D_\frac{\eta}{\sqrt{\lambda}}$, where 
$\eta = \ln(1+\sqrt{2})$; i.e., where $\sinh \eta = 1$, then
the inequality~(\ref{h1}) holds.  
More generally, for all $K \in {\cal K}(\HH^2_{\lambda})$,
\begin{equation}
P_K^2 \geq A_K \left( 4\pi + \lambda A_K \right),
\label{isoper-h}
\end{equation}
where equality holds if and only if $K$ is a hyperbolic disc.

In analogy to the spherical case, if we let $\kappa \rightarrow 0^-$,
so that $\lambda \rightarrow 0^+$, then 
the inequalities of 
Theorem~\ref{omnibus-h} yield the classical Bonnesen-type inequality~(\ref{bonneu})
for the Euclidean plane, while~(\ref{isoper-h}) reduces to the classical isoperimetric
inequality~(\ref{isoper}).

{\bf Remark:} The kinematic approach to isoperimetric inequalities also leads to
generalizations of~(\ref{bonneu}) for the mixed area $A(K, L)$ of compact convex sets
in $\RR^2$ \cite{Santa,red}.  This mixed area arises in the computation
of the area of the Minkowski sum $K + L$, which is itself a convolution integral 
of functions with respect to the translative group
for $\RR^2$.  Although the surfaces $\SS^2$ and $\HH^2$ do not admit a 
subgroup of isometries analogous to the translations of $\RR^2$,
it may prove worthwhile to consider convolutions over other
subgroups of isometries, leading to associated kinematic formulas,
containment theorems, and Bonnesen-type isoperimetric inequalities.
Bonnesen-type inequalities in higher dimensions remain
elusive, but perhaps recent
generalizations of Hadwiger's containment theorem to dimensions
greater than 2, such as those of Zhou \cite{Zhou1, Zhou2},
may be
helpful in developing discriminant inequalities in higher dimension
similar to those presented in this article. 


\end{document}